\documentclass[12pt]{article}
\usepackage{amssymb}
\usepackage{amsmath}
\usepackage{enumerate}
\usepackage[latin1]{inputenc}
\def\s{\mathbb{S}}
\def\h{\mathbb{H}}
\def\r{\mathbb{R}}
\def\z{\mathbb{Z}}
\def\q{\mathbb{Q}}
\def\c{\mathbb{C}}
\def\p{\mathbb{P}}
\newtheorem{definition}{Definition}
\newtheorem{remark}{Remark}
\newtheorem{theorem}{Theorem}
\newtheorem{proposition}{Proposition}
\newtheorem{corollary}{Corollary}
\newtheorem{lemma}{Lemma}
\begin{document}

\title{Hamiltonian-minimal Lagrangian submanifolds in complex space forms}

\author{Ildefonso Castro\thanks{Research partially supported by a MEC-Feder grant MTM2004-00109.
**Partially supported by  the grant  No. 10131020 of NSFC.} $\,\,
$ \ Haizhong Li$^{**}$ \and Francisco Urbano$^*$}
\date{}
\maketitle
\pagestyle{plain}

\begin{abstract}
Using Legendrian immersions and, in particular, Legendre curves in odd dimensional spheres and anti De Sitter spaces,
we provide a method of construction of
new examples of Hamiltonian-minimal Lagrangian submanifolds in complex projective and hyperbolic spaces,
including explicit one parameter families of embeddings of quotients of certain product manifolds.
In addition, new examples of minimal Lagrangian submanifolds in complex projective and hyperbolic spaces also appear.
Making use of all of them, we get Hamiltonian-minimal and special Lagrangian cones in complex Euclidean space too.
\end{abstract}

\section{Introduction}
Let $(\widetilde{M}^n,J,\langle,\rangle)$ be a Kaehler manifold of complex dimension $n$, where $J$ is the complex structure and
$\langle,\rangle$ the Kaehler metric. The Kaehler 2-form is defined by $\omega(.,.) =\langle J.,. \rangle $.
An immersion $\psi: M^n \rightarrow \widetilde{M}^n$ of an $n$-dimensional manifold $M$ is called {\em Lagrangian} if $\psi^* \omega \equiv 0$. For this type of immersions, $J$ defines a bundle isomorphism between the tangent bundle $TM$ and the normal bundle $T^\perp M$.

A vector field $X$ on $\widetilde{M}$ is a Hamiltonian vector
field if ${\cal L}_X \omega = g \omega$, for some funcion $g\in
C^\infty (\widetilde{M})$, where $\cal L$ is the Lie derivative in
$\widetilde M$. This means that there exists a smooth function
$F:\widetilde{M} \rightarrow \r$ such that
$X=J\widetilde{\nabla}F$, where $\widetilde{\nabla}$ is the
gradient in $\widetilde{M}$. The diffeomorphisms of the flux $\{
\varphi_t \}$ of $X$ satisfy that $\varphi_t^* \omega =
e^{h_t}\omega $, and so they transform Lagrangian submanifolds
into Lagrangian ones.

In this setting, Oh studied in [O] the following natural variational problem.
A normal vector field $\xi$ to a Lagrangian immersion $\psi: M^n \rightarrow \widetilde{M}^n$ is called {\em Hamiltonian} if
$\xi= J\nabla f $
where $f\in C^\infty (M)$ and $\nabla f$ is the gradient of $f$ with respect to the induced metric.
If $f\in  C_0^\infty (M)$ and $\{ \psi_t :M\rightarrow\widetilde{M} \}$ is a variation of $\psi $ with $\psi_0=\psi $ and
$\frac{d}{dt} \left|_{_{t=0}}  \, \psi_t \right. = \xi$,
then the first variation of the volume functional is given by (see [O]):
\[   \frac{d}{dt}_{\left|_{t=0} \right.}  \, {\rm vol} (M,\psi_t^* \langle, \rangle )  = -\int_M  f \, {\rm div \, }JH  \, dM ,    \]
where $H$ is the mean curvature vector of the immersion $\psi$ and
div denotes the divergence operator on $M$. Oh called the critical
points of this variational problem {\em Hamiltonian minimal} (or
{\em H-minimal} briefly) Lagrangian submanifolds, which are
characterized by the third order differential equation div$JH
\!=\!0$. In particular, minimal Lagrangian submanifolds (i.e. with
vanishing mean curvature vector) and, more generally, Lagrangian
submanifolds with parallel mean curvature vector are trivially
H-minimal ones.

Even when $\widetilde{M}$ is a simply-connected complex space form, only few examples of H-minimal Lagrangian submanifolds are known outside the class of Lagrangian submanifolds with parallel mean curvature vector.

This can be a brief history of them: In 1998 it was classified in
[CU1] the $\s^1$-invariant H-minimal Lagrangian tori in complex
Euclidean plane $\c^2$. H-minimal Lagrangian cones in $\c^2$ were
studied in 1999, see [SW]. In 2000 and 2002 (see [HR1] and [HR2])
a Weierstrass type representation formula is derived to describe
all H-minimal Lagrangian tori and Klein bottles in $\c^2$. When
the ambient space is the complex projective plane $\c\p^2$ or the
complex hyperbolic plane $\c\h^2$, conformal parametrizations of
H-minimal Lagrangian surfaces using holomorphic data were obtained
in [HR3] and [HR4] in 2002 and 2003. Making use of this technique,
in 2003 H-minimal Lagrangian symply periodic cylinders and
H-minimal Lagrangian surfaces with a non conical singularity in
$\c^2$ were constructed in [A]. Finally, only very recently we can
find in [M] some examples of H-minimal Lagrangian submanifolds of
arbitrary dimension in $\c^n$ and $\c\p^n$ and in [ACR] a
classification of H-minimal Lagrangian submanifolds foliated by
$(n\!-\!1)$-spheres in $\c^n$ is given.

Our aim in this paper is the construction of H-minimal Lagrangian
submanifolds in complex Euclidean space $\c^n$, the complex
projective space $\c\p^n$ and the complex hyperbolic space
$\c\h^n$, for arbitrary $n\geq 2$. The examples in $\c\p^n$ are
constructed by projections, via the Hopf fibration $\Pi:\s^{2n+1}
\rightarrow \c\p^n$, of certain family of Legendrian submanifolds of
the sphere $\s^{2n+1}$ (Corollary 1). The cones with links in this
family of Legendrian submanifolds provide new examples of
H-minimal Lagrangian submanifolds in $\c^{n+1}$ (Corollary 6).
Using the Hopf fibration $\Pi:\h_1^{2n+1} \rightarrow \c\h^n$ and
a similar family of Legendrian submanifolds of the anti De Sitter
space $\h_1^{2n+1}$ (see Corollary 7), we also find out the
examples of H-minimal Lagrangian submanifolds in $\c\h^n$.

In $\c\p^n$, we emphasize two different one-parameter families of H-minimal Lagrangian immersions described in Corollaries 3 and 5; as a particular case of one of them, in Corollary 4 we provide explicit Lagrangian H-minimal embeddings of certain quotients of $\s^1\times \s^{n_1} \times \s^{n_2}$, $n_1+n_2+1 = n$.

In $\c\h^n$, we also point up in Corollary 8 a one-parameter family of H-minimal Lagrangian immersions, which (in the eaisest cases) induce explicit Lagrangian H-minimal embeddings of certain quotients of $\s^1\times \s^{n_1} \times \r\h^{n_2}$, $n_1+n_2+1 = n$ (see Corollary 9), where $\r\h^{n_2}$ denotes the real hyperbolic space.

As a byproduct, using our method of construction, we also obtain
new examples of minimal Lagrangian submanifolds in $\c\p^n$ (see
Corollary 3, Remark 2 and Corollary 5) and $\c\h^n$ (see
Corollaries 7 and 10), as soon as special Lagrangian cones in
$\c^{n+1}$ (see Corollary 6).

\section{Lagrangian submanifolds versus Legendrian submanifolds}

Let $\c^{n+1}$ be the complex Euclidean space endowed with the Euclidean metric $\langle, \rangle$ and the complex
structure $J$. The Liouville 1-form is given by $\Lambda_z (v) = \langle v, J z \rangle $,
 $\forall z \in \c^{n+1}, \forall v \in T_z \c^{n+1}$, and the Kaehler 2-form is $\omega = d\Lambda /2 $.
 We denote the $(2n\!+\!1)$-dimensional unit sphere in $\c^{n+1}$ by $\s^{2n+1}$ and
 by $\Pi:\s^{2n+1}\rightarrow\c\p^n$, $\Pi(z)=[z]$, the Hopf fibration of $\s^{2n+1}$ on
 the complex projective space $\c\p ^{n}$. We also denote the Fubini-Study metric, the complex structure and the K\"{a}hler two-form in $\c\p^n$ by $\langle, \rangle$, $J$ and $\omega$ respectively. This metric has constant holomorphic sectional curvature 4.

We will also note by $\Lambda $ the restriction to $\s^{2n+1}$ of
the Liouville 1-form of $\c^{n+1}$. So $\Lambda $ is the contact
1-form of the canonical Sasakian structure on the sphere
$\s^{2n+1}$. An immersion $\phi :M^n \rightarrow \s^{2n+1}$ of an
$n$-dimensional manifold $M$ is said to be {\em Legendrian} if
$\phi^* \Lambda \equiv 0$. So $\phi $ is isotropic in $\c^{n+1}$,
i.e. $\phi^* \omega \equiv 0$ and, in particular, the normal
bundle $T^\perp M = J (TM) \oplus {\rm span \, } \{ J\phi \} $.
This means that $\phi $ is horizontal with respect to the Hopf
fibration $\Pi:\s^{2n+1} \rightarrow \c\p^n$ and, hence $\Phi= \Pi
\circ \phi: M^n\rightarrow \c\p^n$ is a Lagrangian immersion and
the induced metrics on $M^n$ by $\phi$ and $\Phi $ are the same.
It is easy to check that $J\phi $ is a totally geodesic normal
vector field and so the second fundamental forms of $\phi $ and
$\Phi $ are related by
\[
\Pi _* (\sigma_\phi (v,w)) = \sigma_\Phi (\Pi_* v, \Pi_* w) , \forall v,w \in TM .
\]
So the mean curvature vector $H$ of $\phi $ satisfies that $\langle H, J \phi \rangle =0$  and,
in particular, {\em $\phi :M^n \rightarrow \s^{2n+1}$ is minimal if and only if $\Phi= \Pi \circ \phi: M^n\rightarrow \c\p^n $ is minimal.}

In this way, we can construct (minimal) Lagrangian submanifolds in $\c\p^n$ by projecting Legendrian ones in $\s^{2n+1}$ by the Hopf fibration $\Pi$.

Conversely, it is well known that if $\Phi:M^n\longrightarrow \c\p^{n}$  is a Lagrangian immersion, then $\Phi$ has a horizontal {\em local} lift to $\s^{2n+1}$ with respect to the Hopf fibration $\Pi$, which is unique up to rotations.
We note that only Lagrangian immersions in $\c\p^n$ have this type of lifts.

In this article, we will construct examples of Lagrangian
submanifolds of $\c\p^n$ by constructing examples of Legendrian
submanifolds of $\s^{2n+1}$. Thus we study now some geometric
properties of the Legendrian submanifolds in $\s^{2n+1}$.
\vspace{0.2cm}

Let $\Omega $ be the complex $n$-form on $\s^{2n+1}$ given by
\[
\Omega_z(v_1,\dots,v_n)=\det_{\,\c} \, \{ z,v_1,\dots,v_n \} .
\]
If $\phi :M^n\rightarrow \s^{2n+1}$ is a Legendrian immersion of a manifold $M$,
then $\phi^* \Omega $ is a complex $n$-form on $M$. In the following result we analyze this $n$-form $\phi^* \Omega$.
\begin{lemma}
If $\phi :M^n\rightarrow \s^{2n+1}$ is a Legendrian immersion of a manifold $M$, then
\begin{equation} \label{eq:covariante}
\nabla (\phi^* \Omega)= \alpha_H \otimes \phi^* \Omega,
\end{equation}
where $\alpha_H $ is the one-form on $M$ defined by $\alpha_H (v)= n \,i \langle H, Jv \rangle$
and $H$ is the mean curvature vector of $\phi$.
Consequently, if $\phi $ is minimal then $M$ is orientable.
\end{lemma}
{\it Proof:\/} Let $\{ E_1, \dots, E_n \}$ be an orthonormal frame
on an open subset $U\subset M$, $p\in U$ such that $\nabla_v E_i
=0$, $\forall v\in T_p M$, $i=1,\dots ,n$. We define $A:U
\rightarrow U(n+1)$ by $A =\{ \phi,\phi_* (E_1), \dots , \phi_*
(E_n)  \}$. Then
\[
(\nabla_v \phi^* \Omega)(E_1,\dots,E_n)=v (\det_\c A)= \det_\c A \, \, {\rm Trace \, } (v(A) \bar{A}^t),
\]
where $\bar{A}^t$ denotes the transpose conjugate matrix of $A$. We easily have that
\[
v(A)=\{\phi_{*}(v),\sigma_\phi(v,E_1(p))-\langle v,E_1(p)\rangle\phi,\dots,\sigma_\phi(v,E_{n}(p))-\langle v,E_{n}(p)\rangle\phi\},
\]
and so we deduce that
\[ (\nabla_v \phi^* \Omega)(E_1(p),\dots,E_n(p)) = n \, i \langle H(p) ,Jv\rangle (\phi^* \Omega)(E_1,\dots,E_n)(p). \]
Using this in the above expression we get the result.$_\diamondsuit$
\vspace{0.2cm}

Suppose that our Legendrian submanifold $M$ is oriented. Then we can consider the well defined map given by
\[
\begin{array}{c}
\beta : M^n \longrightarrow \r / 2\pi \z
\\
e^{i\beta (p)}=(\phi^* \Omega)_p (e_1,\dots,e_n)
\end{array}
\]
where $\{ e_1,\dots, e_n \}$ is an oriented orthonormal frame in
$T_pM$. We will call $\beta$ the {\em Legendrian angle} map of
$\phi $. As a consequence of  (\ref{eq:covariante}) we obtain
\begin{equation} \label{eq:gradiente}
J\nabla \beta = n H,
\end{equation}
and so we deduce the following result.
\begin{proposition}
Let $\phi :M^n\rightarrow \s^{2n+1}$ be a Legendrian immersion of an oriented manifold $M$. Then $\phi $ is minimal if and only if the Legendrian angle map $\beta $ of $\phi $ is constant.
\end{proposition}
\vspace{0.3cm}

On the other hand, a vector field $X$ on $\s^{2n+1}$ is a {\em contact vector field} if ${\cal L}_X \Lambda = g \Lambda $, for some function
$g\in C^\infty (\s^{2n+1})$, where $\cal L$ is the Lie derivative in $\s^{2n+1}$. It is well known (see [McDS]) that $X$ is a contact vector field if and only if there exists $F\in C^\infty(\s^{2n+1})$ such that
\[
X_z=J(\overline{\nabla}F)_z + 2F Jz, \quad z\in \s^{2n+1} ,
\]
where $\overline{\nabla}F$ is the gradient of $F$.
The diffeomorphisms of the flux $\{ \varphi_t \}$ of $X$ are contactmorphisms of $\s^{2n+1}$, that is,
 $\varphi_t^* \Lambda = e^{h_t}\Lambda $, and so they transform Legendrian submanifolds into Legendrian ones. The Lie algebra of the group of contactmorphisms of $\s^{2n+1}$ is the space of contact vector fields. In this setting, it is natural to study the following variational problem.

Let $\phi :M^n\rightarrow \s^{2n+1}$ a Legendrian immersion with mean curvature vector $H$. A normal vector field $\xi_f$ to $\phi $ is called  a {\em contact field} if
\[
\xi_f = J\nabla f + 2 f J\phi,
\]
where $f\in C^\infty (M)$ and $\nabla f$ is the gradient of $f$
with respect to the induced metric. If $f\in  C_0^\infty (M)$ and
$\{ \phi_t :M\rightarrow\s^{2n+1} \}$ is a variation of $\phi $
with $\phi_0=\phi $ and $\frac{d}{dt}_{\left|_{t=0}\right.} \phi_t
= \xi_f $, the first variation of the volume functional is given
by
\[
\frac{d}{dt}_{\left|_{t=0} \right.}  \, {\rm vol} (M,\phi_t^* \langle, \rangle )  = -\int_M \langle  H, \xi_f \rangle \, dM .    \]
But using the Stoke's Theorem,
\[
\begin{array}{c}
\int_M \langle  H, \xi_f \rangle \, dM = \int_M \langle  H, J\nabla f + 2 f J\phi \rangle \, dM
\\ \\
=-\int_M \langle  JH, \nabla f  \rangle \, dM = \int_M  f {\rm div} \, JH \, dM .
\end{array}
\]
This means that the critical points of the above variational problem are Legendrian submanifolds such that
\[ {\rm div}  JH =0. \]
We name the critical points of this variational problem in the following definition.
\begin{definition}
A Legendrian immersion $\phi :M^n \rightarrow \s^{2n+1}$ is said to be {\em  contact minimal} (or briefly {\em C-minimal}) if $ {\rm div}  JH =0$.
\end{definition}
Clearly, minimal Legendrian submanifolds and Legendrian submanifolds with parallel mean curvature vector are C-minimal.
As a consequence of (\ref{eq:gradiente}) and the geometric relationship between Legendrian and Lagrangian submanifolds mentioned at the beginning of this section, we get the following.
\begin{proposition}
If $\phi :M^n\rightarrow \s^{2n+1}$ is a Legendrian immersion of a Riemannian manifold $M$,
then:
\begin{enumerate}
\item If $M$ is oriented, $\phi $ is C-minimal if and only if the Legendrian angle $\beta $ of $\phi $ is a harmonic map.
\item $\phi $ is C-minimal if and only if $\Phi= \Pi \circ \phi :M^n\rightarrow \c\p^n$ is H-minimal.
\end{enumerate}
\end{proposition}

\section{A new construction of C-minimal Legendrian immersions}

After Proposition 2, it is clear that constructing C-minimal Legendrian immersions in odd-dimensional spheres is a good way to find out H-minimal Lagrangian submanifolds in $\c\p^n$. This is the purpose of this section. But first we need to introduce some notations.
Let $n_1$ and $n_2$ be nonnegative integer numbers and $n=n_1+n_2+1$.
If $SO(m)$ denotes the special orthogonal group, then
$SO(n_1+1)\times SO(n_2+1)$ acts on $\s^{2n+1}\subset \c^{n+1}$, $n=n_1+n_2+1$, as a subgroup of isometries in the following way:
\begin{equation}
\label{eq:action1}
(A_1,A_2)\in SO(n_1+1)\times SO(n_2+1)\longmapsto \left(\begin{array}{c|c}
A_1 &  \\
\hline
&  A_2
\end{array}\right) \in SO(n+1).
\end{equation}
\begin{theorem}\
Let $\psi_i : N_i \rightarrow \s^{2n_i+1}\subset \c^{n_i+1}$ be
Legendrian isometric immersions of $n_i$-dimensional oriented
Riemannian manifolds $( N_i, g_{N_i} )$, $i=1,2$, and
$\gamma=(\gamma_1,\gamma_2):I\rightarrow \s^3\subset \c^2$  be a
Legendre curve. Then the map
\[
\phi:I\times N_1 \times N_2 \longrightarrow \s^{2n+1}\subset \c^{n+1}=\c^{n_1+1}\times \c^{n_2+1}, \, n=1+n_1+n_2 ,
\]
 defined by
\[
\phi(s,p,q)=\left(\gamma_1(s)\psi_1(p),\gamma_2(s)\psi_2(q)\right)
\]
is a Legendrian immersion in $\s^{2n+1}$ whose induced metric is
\begin{equation}
\label{eq:metric1} \langle , \rangle = |\gamma'|^2 ds^2 +
|\gamma_1|^2 g_{N_1} + |\gamma_2|^2 g_{N_2}
\end{equation}
and whose Legendrian angle map is
\begin{equation}
\label{eq:angle1}
\beta_\phi \equiv n_1 \pi +\beta_\gamma + n_1 \arg \gamma_1 + n_2 \arg \gamma_2 + \beta_{\psi_1}
+\beta_{\psi_2} \quad \mod \, 2\pi,
\end{equation}
where $\beta_\gamma$ denotes the Legendre angle of $\gamma$ and $\beta_{\psi_i}$ the Legendrian angle map of $\psi_i$, $i=1,2$.

Moreover, a Legendrian immersion $\phi :M^n\longrightarrow \s^{2n+1}$  is  invariant under the action (\ref{eq:action1}) of $SO(n_1\!+\!1)\times SO(n_2\!+\!1)$, with $n=n_1+n_2+1$ and $n_1,n_2 \geq 2$, if and only if
$\phi $ is locally congruent to one of the above Legendrian immersions when $\psi_i$ are the totally geodesic Legendrian embeddings of $\s^{n_i}$ in $\s^{2n_i+1}$, $i=1,2$; that is,
$\phi $ is locally given by $\phi(
s,x,y)=(\gamma_1(s)\, x,\gamma_2(s)\, y)$, $x\in \s^{n_1}$, $y\in \s^{n_2}$, for a certain Legendre curve $\gamma$ in $\s^3$.
\end{theorem}
These Legendrian immersions introduced in Theorem 1 have singularities in the points $(s,p,q)\in I \times N_1 \times N_2$ where either $\gamma_1(s)=0$ or $\gamma_2(s)=0$.
\vspace{0.2cm}

{\it Proof:\/}
If $'$ denotes derivative with respect to $s$, and $v$ and $w$ are arbitrary tangent vectors to $N_1$ and $N_2$ respectively, it is clear that
\begin{eqnarray}\label{eq:derivadas}
\phi_s=\phi_* (\partial_s,0,0)=(\gamma_1' \, \psi_1, \gamma_2'\, \psi_2), \nonumber \\
\phi_* (v) := \phi_* (0,v,0)=(\gamma_1 \, \psi_{1_*} (v), 0), \\
\phi_* (w) := \phi_* (0,0,w) = (0,\gamma_2 \, \psi_{2_*} (w)). \nonumber
\end{eqnarray}
Recall that $g_{N_1}$ and  $g_{N_2}$ are the induced metrics on
$N_1$ and $N_2$ by $\psi_1$ and  $\psi_2$ respectively. From
(\ref{eq:derivadas}) and using that $\psi_1$ and $\psi_2$ are
Legendrian immersions, the induced metric on $I\times N_1\times
N_2$ by $\phi$ is given by $|\gamma'|^2 ds^2 + |\gamma_1|^2
g_{N_1} + |\gamma_2|^2 g_{N_2}$.
From the Legendrian characters of $\gamma$, $\psi_1$ and $\psi_2$,
it follows that the immersion $\phi$ is also Legendrian.

In order to compute the Legendrian angle map $\beta_{\phi}$, let
$\{e_1,\dots,e_{n_1}\}$ and $\{e'_1,\dots,e'_{n_2}\}$  be oriented
local orthonormal frames on $N_1$ and $N_2$ respectively. Then
\begin{equation}
\label{eq:bo}
\{u_1,v_1,\dots,v_{n_1},w_1,\dots,w_{n_2} \}
\end{equation}
 defined by
\begin{eqnarray*}
u_1&=&\left(\frac{\partial_s}{|\gamma'|},0,0\right) \\
v_j&=&\left(0,\frac{e_{j}}{|\gamma_1|},0\right), \, 1\leq j \leq n_1, \\
w_k&=&\left(0,0,\frac{e'_{k}}{|\gamma_2|}\right),\, 1 \leq k \leq n_2 ,
\end{eqnarray*}
is an oriented local orthonormal frame on $I \times N_1 \times
N_2$. Putting $\phi=\gamma_1(\psi_1,0)+\gamma_2(0,\psi_2)$ and
$\phi_*(u_1)=\frac{\gamma_1'}{|\gamma'|}(\psi_1,0)+\frac{\gamma_2'}{|\gamma'|}(0,\psi_2)$,
we have that
\[
\begin{array}{c}
e^{i\beta_{\phi}}=\det_\c \, \{ \phi, \phi_*(u_1),\dots,\phi_*(v_j),\dots, \phi_*(w_k), \dots \}= \\
\frac{\gamma_1^{n_1}\gamma_2^{n_2}(\gamma_1 \gamma_2'
-\gamma_1' \gamma_2)}{|\gamma'||\gamma_1|^{n_1}|\gamma_2|^{n_2}} 
\det_\c \, \{ (\psi_1,0),(0,\psi_2), \dots,
(\psi_{1_*}(e_j),0),\dots,(0,\psi_{2_*}(e'_{k}), \dots \}  .
\end{array}
\]
In this way we obtain that
\[
e^{i\beta_{\phi}(s,p,q)}=(-1)^{n_1}\, e^{i (n_1\arg\gamma_1+n_2\arg\gamma_2)(s)}\,
\frac{(\gamma_1 \gamma_2'-\gamma_1' \gamma_2)(s)}{|\gamma'(s)|}\,
\det_{\c} A_1(p)\,\det_{\c} A_2(q),
\]
where $A_1$ and $A_2$ are the matrices
\[
A_1=\{\psi_1, \psi_{1_*}(e_1),\dots,\psi_{1_*}(e_{n_1})\}
\]
and
\[
A_2=\{\psi_2,\psi_{2_*}(e'_1),\dots,\psi_{2_*}(e'_{n_2})\}.
\]
Taking into account the definition of the Legendrian angle map given in section 2, we finally arrive at
\[
e^{i\beta_{\phi}(s,p,q)}=(-1)^{n_1}\, e^{i (\beta_\gamma +n_1\arg\gamma_1+n_2\arg\gamma_2)(s)}\,
e^{i\beta_{\psi_1}(p)} \, e^{i\beta_{\psi_2}(q)}. \]
This proves the first part of the result.
\vspace{0.5cm}

On the other hand, let
$\phi:M^n\rightarrow\s^{2n+1}\subset\c^{n+1}$ be a Legendrian
immersion  which is invariant under the action (\ref{eq:action1})
of $SO(n_1+1)\times SO(n_2+1)$, $n=n_1+n_2+1$. Let $p$ be any
point of $M$ and let $z=(z_1,\dots,z_{n+1})=\phi(p)$. As $\phi$ is
invariant under the action of $SO(n_1+1)\times SO(n_2+1)$, for any
matrix $X=(X_1,X_2)$ in the Lie algebra of $SO(n_1+1)\times
SO(n_2+1)$, the curve $t\mapsto ze^{t\hat X}$ with
\[
\hat X=\left( \begin{array}{c|c}
\mbox{X}_1 &  \\
\hline
& \mbox{X}_2
\end{array}\right)
\]
lies in the submanifold. Thus its tangent vector at $t=0$ satisfies
\[
z\hat X\in\phi_*(T_pM).
\]
Since $\phi$ is a Legendrian immersion, this implies that
\[
\Im (z\hat X\hat Y\bar{z}^t)=0
\]
for any matrices $X=(X_1,X_2)$, $Y=(Y_1,Y_2)$ in the Lie algebra of $SO(n_1+1)\times SO(n_2+1)$.
As $n_1+1\geq 3$ and $n_2+1\geq 3$,  it is easy to see from the last equation that $\Re (z_1,\dots,z_{n_1+1})$ and $\Im (z_{1},\dots,z_{n_1+1})$ (respectively $\Re (z_{n_1+2},\dots,z_{n+1})$ and $\Im (z_{n_1+2},\dots,z_{n+1})$) are linear dependent. As $SO(n_1+1)$ acts transitively on $\s^{n_1}$ and $SO(n_2+1)$ acts transitively on $\s^{n_2}$, we obtain that $z$ is in the orbit (under the action of $SO(n_1+1)\times SO(n_2+1)$ described above) of the point $(z^0_1,0,\dots,0,z^0_{n_1+2},0,\dots,0)$, with $|z^0_1|^2=\sum_{i=1}^{n_1+1}|z_i|^2$ and $|z^0_{n_1+2}|^2=\sum_{j=n_1+2}^{n+1}|z_j|^2$. This implies that locally $\phi$ is the orbit under the action of $SO(n_1+1)\times SO(n_2+1)$ of a curve $\gamma$ in $\c^2\equiv \c^n\cap\{z_2=\dots=z_{n_1+1}=z_{n_1+3}=\dots=z_{n+1}=0\}$. Therefore $M$ is locally $I\times\s^{n_1}\times\s^{n_2}$, with $I$ an interval in $\r$. Moreover, $\phi$ is given by
\[
\phi(s,x,y)=(\gamma_1(s)\, x,\gamma_2(s) \, y),
\]
where $\gamma=(\gamma_1,\gamma_2)$ must be a Legendre curve in
$\s^3\subset\c^2$. Finally, as $\phi$ is a Legendrian submanifold,
the result follows using the first part of this
Theorem.$_{\diamondsuit}$ \vspace{0.3cm}

In the following result we make use of the method described in
Theorem 1 to obtain new minimal and C-minimal Legendrian
immersions, which will provide (projecting via the Hopf fibration)
new non-trivial minimal and H-minimal immersions in $\c\p^n$.

\begin{corollary}\
Let $\psi_i : N_i \longrightarrow \s^{2n_i+1}$, $i=1,2$, be C-minimal Legendrian immersions
of $n_i$-dimensional oriented Riemannian manifolds $N_i$, $i=1,2$, and $\gamma=(\gamma_1,\gamma_2):I\rightarrow \s^3\subset \c^2$  be a Legendre curve.
Then the Legendrian immersion described in Theorem 1 given by
\begin{eqnarray*}
\phi : I \times N_1\times N_2 \longrightarrow \s^{2n+1}, \, n=n_1+n_2+1,  \\
\phi(t,p,q)=\left(\gamma_1(t)\psi_1(p),\gamma_2(t)\psi_2(q)\right)
\end{eqnarray*}
is  C-minimal  if and only if $(\gamma_1,\gamma_2)$ is a solution of some equation in the two parameter family of o.d.e.
\begin{equation}
\label{eq:gamma12}
(\gamma_1' \overline{\gamma_1})(t) = -( \gamma_2' \overline{\gamma_2})(t)=
 - \, e^{i(\lambda+\mu t)} \,    \overline{\gamma_1}(t)^{n_1+1} \,
\overline{\gamma_2}(t)^{n_2+1}, \, \lambda, \mu \in \r .
\end{equation}
Moreover, the above Legendrian immersion $\phi $ is minimal if and only if
$\psi_i$, $i=1,2$, are  minimal and $(\gamma_1,\gamma_2)$ is a solution of some o.d.e. of
(\ref{eq:gamma12}) with $\mu =0$.
\end{corollary}
\begin{remark}
{\rm If we make a $\theta$-rotation of a Legendre curve $\gamma$
solution of (\ref{eq:gamma12}) for the parameters $(\lambda,\mu)$,
the new Legendre curve is a solution  of (\ref{eq:gamma12}) for
the parameters $(\lambda -(n+1)\theta,\mu)$. The corresponding
immersions given in Corollary 1 are related by
$\tilde{\phi}=e^{i\theta}\phi$ and so they are congruent. In this
way, taking $\theta=\frac{\pi/2+\lambda}{n+1}$, up to congruences
it is sufficient to consider solutions of the one parameter family
of equations
\begin{equation}
\label{eq:gammamu}
(\gamma_j' \overline{\gamma_j})(t) =
 (-1)^{j-1} i\, e^{i\mu t} \,    \overline{\gamma_1}(t)^{n_1+1} \,
\overline{\gamma_2}(t)^{n_2+1}, \,  \mu \in \r , \, j=1,2.
\end{equation}
}
\end{remark}
{\it Proof:\/} Recall from (\ref{eq:angle1}) that
\[
\beta_\phi \equiv n_1 \pi +\beta_\gamma + n_1 \arg \gamma_1 + n_2 \arg \gamma_2 + \beta_{\psi_1}
+\beta_{\psi_2} \quad \mod \, 2\pi,
\]
where
$\phi $ is one of the Legendrian immersions described in Theorem 1.

Using Proposition 2, $\phi $ is C-minimal if and only if
$\Delta\beta_\phi =0$. So we must compute the Laplacian of
$\beta_\phi$.  For this purpose we use the orthonormal frame
(\ref{eq:bo}) and after a long but direct computation we obtain
\begin{equation}
\label{eq:Laplace}
\Delta\beta_\phi = \frac{1}{|\gamma'|^2}
\left(
\frac{\partial^2\beta_\phi}{\partial s^2}
+\frac{{\rm d}}{{\rm d}s}\left( \log \frac{|\gamma_1|^{n_1} |\gamma_2|^{n_2}}{|\gamma'|} \right) \frac{\partial\beta_\phi}{\partial s}
\right)
+\frac{\Delta_1\beta_{\psi_1}}{|\gamma_1|^2}
+\frac{\Delta_2\beta_{\psi_2}}{|\gamma_2|^2},
\end{equation}
where $\Delta_i$ are the Laplace operators in $( N_i, g_{N_i} )$, $i=1,2$.

The assumptions of the Corollary 1 imply that
$\Delta_i\beta_{\psi_i}=0$, $i=1,2$, using Proposition 2 again.
So $\phi$ is C-minimal if and only if
\begin{equation}\label{eq:Hminbeta}
\frac{\partial^2\beta_\phi}{\partial s^2}
+\frac{{\rm d}}{{\rm d}s}\left( \log \frac{|\gamma_1|^{n_1} |\gamma_2|^{n_2}}{|\gamma'|} \right) \frac{\partial\beta_\phi}{\partial s} =0.
\end{equation}
From (\ref{eq:metric1}), we have that $\gamma_i(0)\neq 0$, $i=1,2$, since we want $\phi $ to be regular.
So we can choose, up to  reparametrizations, $\gamma=\gamma(t)$ to satisfy that
$|\gamma'(t)|=|\gamma_1(t)|^{n_1} |\gamma_2(t)|^{n_2} $.
Thus (\ref{eq:Hminbeta}) becomes
\[
\frac{\partial^2\beta_\phi}{\partial t^2} =0.
\]
This means that $\beta_\phi (t,p,q)=f(p,q)+t \, g(p,q)$, for certain functions $f,g$ defined on $N_1 \times N_2$.
Using (\ref{eq:angle1}), we obtain that $g(p,q)=$constant and
\begin{equation}
\label{eq:Cminimality}
(\beta_\gamma + n_1 \arg \gamma_1 + n_2 \arg \gamma_2)(t)=\lambda + \mu t, \, \lambda, \mu  \in \r.
\end{equation}
The definition of  the Legendrian angle $\beta_\gamma$ of  $\gamma$ is given, in particular, by
\[
e^{i \beta_\gamma}=\frac{1}{|\gamma'|}(\gamma_1 \gamma_2' -\gamma_2 \gamma_1').
\]
Using this, it is easy to check that
(\ref{eq:Cminimality}) can be written as
\[
\gamma_1' \overline{\gamma_1} = - \gamma_2' \overline{\gamma_2}=
 - \, e^{i(\lambda+\mu t)} \,    \overline{\gamma_1}^{n_1+1} \,
\overline{\gamma_2}^{n_2+1},
\]
that is exactly (\ref{eq:gamma12}).

Finally, using Proposition 1, $\phi $ is minimal if and only if $\beta_\phi$ is constant. This is equivalent to that
$\beta_{\psi_i}$, $i=1,2$, are constant (i.e. $\psi_i$ are minimal from Proposition 1 again) and $\beta_\gamma + n_1 \arg \gamma_1 + n_2 \arg \gamma_2$ is constant. But this corresponds to the case $\mu =0$ in (\ref{eq:Cminimality}) and so to the case $\mu =0$ in (\ref{eq:gamma12}).$_\diamondsuit$
\vspace{0.2cm}

It is rather difficult to describe the general solution of (\ref{eq:gammamu}).
However it is an exercise to check that for any $\delta \in (0,\pi/2)$ the Legendre curve
\begin{equation}\label{eq:gammadelta}
\gamma_\delta (t)=( c_{\delta} \, \exp(i s_{\delta}^{n_1+1}  c_{\delta}^{n_2-1}    t) ,
s_{\delta} \, \exp (-is_{\delta}^{n_1-1} c_{\delta}^{n_2+1}  t)),
\end{equation}
satisfies (\ref{eq:gammamu}) for
$\mu = s_{\delta}^{n_1-1} c_{\delta}^{n_2-1} \left( (n_1+1)s_{\delta}^2  -(n_2+1)c_{\delta}^2  \right)$, where $c_{\delta}=\cos\delta$, $s_{\delta}=\sin\delta$.
We observe that this value of $\mu $ vanishes if and ond only if $\tan^2 \delta = (n_2+1)/(n_1+1)$.
In this way we are able to obtain the following explicit family of examples.

\begin{corollary}\
Let $\psi_i : N_i \longrightarrow \s^{2n_i+1}$, $i=1,2$, be C-minimal Legendrian immersions
of $n_i$-dimensional Riemannian manifolds $N_i$, $i=1,2$.

Given any $\delta \in (0,\pi/2)$ and denoting $c_{\delta}=\cos\delta$ and $ s_{\delta}=\sin\delta$, the map
\[
\begin{array}{c}
\phi_\delta: \r \times N_1\times N_2 \longrightarrow \s^{2n+1}, \,
\, n=n_1+n_2+1,
\\  \\
\phi_\delta (t,p,q)=
( c_{\delta} \, \exp( i s_{\delta}^{n_1+1} c_{\delta}^{n_2-1}   t) \,\psi_1 (p) \,  , \,
 s_{\delta} \, \exp( -i s_{\delta}^{n_1-1} c_{\delta}^{n_2+1}   t) \,\psi_2 (q))
\end{array}
\]
is a C-minimal Legendrian immersion.

In particular, using minimal Legendrian immersions $\psi_i$, $i=1,2$, and $\delta_0 = \arctan \sqrt{(n_2+1)/(n_1+1)}$,
the Legendrian immersion $\phi_{\delta_0}: \r \times N_1\times N_2 \longrightarrow \s^{2n+1}, \, n=n_1+n_2+1$,
is minimal.
\end{corollary}
{\it Proof:\/} We simply remark that we do not need the orientability assumption because, in this case, it is easy to check that the Legendrian immersions $\phi_\delta$ satisfy div$JH=0$ and so they are C-minimal (see Definition 1).$_\diamondsuit$
\vspace{0.1cm}

To finish this section, we pay now our attention to the equation
(\ref{eq:gammamu}) with $\mu=0$. We observe that it is exactly
equation (6) in [CU2, Lemma 2] (in the notation of that paper, put
$p=n_1$ and $q=n_2$). If we choose the initial conditions
$\gamma(0)=(\cos \theta, \sin \theta)$, $\theta \in (0,\pi/2)$, we
can make use of the study made in [CU2].

\begin{lemma}
Let $\gamma_\theta=(\gamma_1,\gamma_2):I\subset \r \rightarrow \s^3$ be the only curve solution of
\begin{equation}\label{eq:gamma0}
\gamma'_j\bar{\gamma}_j=(-1)^{j-1}\, i\, \bar{\gamma}_1^{n_1+1}\bar{\gamma}_2^{n_2+1}, \, j=1,2,
\end{equation}
satisfying the real initial conditions
$\gamma_\theta(0)=(\cos \theta, \sin \theta)$, $\theta \in (0,\pi/2)$.
Then:
\begin{enumerate}
\item $\Re (\gamma_1^{n_1+1}\gamma_2^{n_2+1})=\cos^{n_1+1}\theta \, \sin^{n_2+1}\theta$.
\item For $j=1,2$, $\bar{\gamma}_j(t)=\gamma_j(-t),\,\forall t\in I$.
\item  The functions $|\gamma_j|,\,j=1,2,$ are periodic with the same period $T=T(\theta)$.
Moreover, $\gamma_\theta$ is a closed curve if and only if
\[
\theta \in 
\left\{ \theta \in (0,\frac{\pi}{2}) \, / \, \frac{\cos^{n_1+1}\theta \sin^{n_2+1}\theta}{2\pi}\left ( \int_0^{T}\frac{dt}{|\gamma_1|^2(t)},\int_0^{T}\frac{dt}{|\gamma_2|^2(t)}\right ) \in\q^2  \right\}.
\]
\item If $\theta=\arctan \sqrt{\frac{n_2+1}{n_1+1}}$,
the curve $\gamma_\theta$ is exactly the curve $\gamma_{\delta_0}$ (see Corollary 2 and (\ref{eq:gammadelta})).
\end{enumerate}
\end{lemma}

{\it Proof:\/} Parts 1 and 2 follow directly from parts 2 and 3 in [CU2, Lemma 2].
To prove 3 we define $f(\theta):=\cos^{2(n_1+1)}\theta \sin^{2(n_2+1)}\theta$, $\theta\in(0,\pi/2)$.
It is easy to prove that $f(\theta)\leq (n_1+1)^{n_1+1}(n_2+1)^{n_2+1}/(n+1)^{n+1}$ and the equality holds if and only if $\theta = \delta_0$.
Using this in
parts 4 and 5 in [CU2, Lemma 2], we finish the proof.$_\diamondsuit$




\section{H-minimal Lagrangian submanifolds in complex projective space}

In section 2 we explained that we can construct (minimal, H-minimal) Lagrangian submanifolds in $\c\p^n$
by projecting (minimal, C-minimal) Legendrian submanifolds in $\s^{2n+1}$ by the Hopf fibration $\Pi: \s^{2n+1} \rightarrow \c\p^n$ (see Proposition 2). The aim of this section is to analyze the Lagrangian immersions in $\c\p^n$ that we obtain just by projecting the Legendrian ones deduced in section 3.

First we mention that
if $n_2=0$ in Theorem 1, projecting by the Hopf fibration $\Pi$
we obtain the Examples 1 given in [CMU1].
In this sense, the construction given in Theorem 1 can be considered as generalization of the family introduced in [CMU1].
Some applications of our construction of Theorem 1 when $n=3$ have been used very recently in [MV] to the study of minimal Lagrangian submanifolds in $\c\p^3$.

The Legendrian immersions described in Corollary 1 provide new examples of Lagrangian H-minimal immersions in $\c\p^n$
when we project them by $\Pi$.
If we consider the particular case $n_2=0\Leftrightarrow n_1=n-1$ in the minimal case of Corollary 1, we recover (projecting via the Hopf fibration $\Pi$) the minimal Lagrangian submanifolds of $\c\p^n$ described in [CMU2, Proposition 6], although we used there an unit speed parametrization for $\gamma$.

We write more in detail what we obtain by this procedure if we consider the special case coming from Corollary 2.

\begin{corollary}
Let $\psi_i : N_i \longrightarrow \s^{2n_i+1}$, $i=1,2$, be C-minimal Legendrian immersions
of $n_i$-dimensional Riemannian manifolds $N_i$, $i=1,2$ and $\delta \in (0,\pi/2)$. Then
\[
\Phi_\delta: \s^1 \times N_1 \times N_2 \longrightarrow \c\p^{n}, \, n=n_1+n_2+1,
\]
 given by
\[
\Phi_{\delta}(e^{is},p,q)=[(\cos\delta\,\exp (i s \sin^2 \delta )
\psi_1(p), \sin\delta\,\exp (-is\cos^2 \delta ) \psi_2(q))]
\]
is a H-minimal Lagrangian immersion.

Moreover, $\Phi_{\delta}$ is  minimal if and only if $\psi_i$, $i=1,2$, are minimal and 
$\tan^2 \delta = (n_2+1)/(n_1+1)$.
\end{corollary}

{\it Proof:\/} We consider the C-minimal Legendrian immersions
\[ \phi_\delta: \r \times N_1\times N_2 \longrightarrow \s^{2n+1} \]
given in Corollary 2. Projecting by the Hopf fibration  $\Pi: \s^{2n+1} \rightarrow \c\p^n$ and using Proposition 2
\[ \Pi \circ \phi_\delta: \r \times N_1\times N_2 \longrightarrow \c\p^{n} \]
is a one parameter family of H-minimal Lagrangian immersions.  We analyse when $\Pi \circ \phi_\delta$ is periodic
in its first variable. It is easy to obtain that there exists $A>0$ such that $(\Pi \circ \phi_\delta)(t+A,p,q)=(\Pi \circ \phi_\delta)(t,p,q)$, $\forall (t,p,q)\in \r \times N_1 \times N_2$
if and only if there exists $\theta \in \r$ verifying
\[ \exp( i s_{\delta}^{n_1+1} c_{\delta}^{n_2-1}  A) = e^{i\theta} = \exp( -i s_{\delta}^{n_1-1} c_{\delta}^{n_2+1}  A). \]
From here we deduce that the smallest period $A$ must be given by $A=2\pi/ (s_{\delta}^{n_1-1} c_{\delta}^{n_2-1})$.
If we define the change of variable
\[
\begin{array}{c}
[0,2\pi] \rightarrow [0,2\pi/ (s_{\delta}^{n_1-1} c_{\delta}^{n_2-1}) ] \\ \\
s \mapsto t= s/ (s_{\delta}^{n_1-1} c_{\delta}^{n_2-1})
\end{array}
\]
the Legendre curve $\gamma_\delta$ given in (\ref{eq:gammadelta}) is written as
\[
\gamma_\delta (s)=( c_{\delta} \, \exp(i s_{\delta}^2 \, s) ,
                    s_{\delta} \, \exp (-i c_{\delta}^2  s) ), \, \, s\in [0,2\pi],
\]
and now it is clear that we arrive at the expression of $\Phi_\delta$.

Taking into account that $ \Pi \circ \phi_\delta$ is minimal if and only if $ \phi_\delta $ is minimal (see Section 2) and
using again Corollary 2, we finish the proof of this result.$_\diamondsuit$

\vspace{0.2cm}

We can even get H-minimal Lagrangian embeddings from
a particular case of Corollary 3.

\begin{corollary}
For each $\delta \in (0,\pi/2)$, the immersion $\Phi_{\delta}$ (given in Corollary 3) where $\psi_i$ are the totally geodesic Legendrian embeddings of \,$\s^{n_i}$ into $\s^{2n_i+1}$, $i=1,2$, provides a H-minimal Lagrangian embedding
\[
\begin{array}{c}
\frac{\textstyle \s^1\times\s^{n_1}\times\s^{n_2}}{\textstyle \z_2\times\z_2} \longrightarrow \c\p^n,\, n=n_1+n_2+1 \\ \\
\overline{(e^{is},x,y)} \longmapsto  [(\cos\delta\,\exp (i s
\sin^2 \delta ) x,\sin\delta\,\exp(-i s \cos^2\delta ) y)]
\end{array}
\]
of the quotient of $\s^1\times\s^{n_1}\times\s^{n_2}$ by the action of the group $\z_2\times\z_2$
where the generators $h_1$ and $h_2$ of $\z_2$ act on $\s^1\times\s^{n_1}\times\s^{n_2}$ in the following way
\[
h_1(e^{is},x,y)=(-e^{is},-x,y),\quad h_2(e^{is},x,y)=(-e^{is},x,-y).
\]
\end{corollary}
{\it Proof:\/} We consider the H-minimal Lagrangian immersions
\[
\begin{array}{c}
\Phi_\delta: \s^1 \times \s^{n_1}\times\s^{n_2} \longrightarrow \c\p^{n}, \, n=n_1+n_2+1,
\\  \\
 \Phi_{\delta}(e^{is},x,y)=[(\cos\delta\,\exp (i \sin^2 \delta \, s)\, x,
                            \sin\delta\,\exp (-i\cos^2 \delta \, s)\, y)].
\end{array}
\]
Let $(e^{is},x,y), (e^{i\hat{s}},\hat{x},\hat{y})\in \s^1 \times \s^{n_1}\times\s^{n_2}$. Then
$\Phi_\delta (e^{is},x,y) = \Phi_\delta (e^{i\hat{s}},\hat{x},\hat{y})$ if and only if $\exists \theta \in \r$ such that
\begin{equation}\label{eq:emb1}
\hat{x} = \exp \left( i (\theta + \sin^2 \delta (s-\hat{s}) ) \right) x, \quad
\hat{y} = \exp \left( i (\theta - \cos^2 \delta (s-\hat{s}) ) \right) y .
\end{equation}
As some coordinate of $x\in \s^{n_1}$ and $y \in \s^{n_2}$ is non null, we deduce that
\begin{equation}\label{eq:emb2}
\epsilon_1:= \exp \left( i (\theta + \sin^2 \delta (s-\hat{s}) ) \right) = \pm 1,
\epsilon_2:= \exp \left( i (\theta - \cos^2 \delta (s-\hat{s}) ) \right) = \pm 1 .
\end{equation}
We distinguish the following cases:

(i) $\epsilon_1 = \epsilon_2=\pm 1$:

From (\ref{eq:emb2}) we get that $e^{i\hat{s}}=e^{is}$
and using (\ref{eq:emb1}) we obtain that $\hat{x}=x$, $\hat{y}=y$ if $\epsilon_1 = \epsilon_2=1$ or
$\hat{x}=-x$, $\hat{y}=-y$ if $\epsilon_1 = \epsilon_2=-1$.
Thus $(e^{i\hat{s}},\hat{x},\hat{y})=(e^{is},x,y)$ or $(e^{i\hat{s}},\hat{x},\hat{y})=(e^{is},-x,-y)= (h_1 \circ h_2) (e^{is},x,y)$.

(ii) $\epsilon_1 = -\epsilon_2 = \pm 1$:

From (\ref{eq:emb2}) we get that $e^{i\hat{s}}=-e^{is}$
and using (\ref{eq:emb1}) we obtain that either
$\hat{x}= x$ and $\hat{y}=-y$ if $\epsilon_1 = -\epsilon_2=1$ and so
$(e^{i\hat{s}},\hat{x},\hat{y})=(-e^{is},x,-y)= h_2 (e^{is},x,y)$
or
$\hat{x}= -x$ and $\hat{y}=y$ if $\epsilon_1 = -\epsilon_2=-1$ and so
$(e^{i\hat{s}},\hat{x},\hat{y})=(-e^{is},-x,y)= h_1 (e^{is},x,y)$.
This reasoning proves the result.$_\diamondsuit$

\begin{remark}
{\rm If $\tan^2 \delta =(n_2+1)/(n_1+1) $ in Corollary 4 we obtain a minimal Lagrangian embedding that
generalizes a well known example $\frac{\s^1 \times \s^{n-1} }{\z_2 } \longrightarrow \c\p^{n}$ studied by H. Naitoh in [N], which corresponds to take $n_2=0\Leftrightarrow n_1=n-1$ in Corollary 4.

We also remark that $h_1$ (resp. $h_2$) preserves the orientation of $\s^1 \times \s^{n_1}\times\s^{n_2}$
if and only if $n_1$ (resp. $n_2$) is odd.
Thus $\frac{\s^1 \times \s^{n_1} \times \s^{n_2}}{\z_2 \times \z_2}$ is an orientable manifold
if and only if $n_1$ and $n_2$ are odd.
}
\end{remark}

We finish this section making use of the information given in Lemma 2 for the solutions of equation (\ref{eq:gammamu}) with $\mu =0$.

Let $\theta \in (0,\pi/2)$ and $\gamma_\theta$ be the only solution of
(\ref{eq:gamma0}) satisfying $\gamma_\theta(0)=(\cos \theta, \sin \theta)$.
We consider the C-minimal Legendrian immersions
\[ \phi_\theta: I \times N_1\times N_2 \longrightarrow \s^{2n+1} \]
constructed with $\gamma_\theta$. Projecting by the Hopf fibration  $\Pi: \s^{2n+1} \rightarrow \c\p^n$ and using Proposition 2
\[ \Pi \circ \phi_\theta: I \times N_1\times N_2 \longrightarrow \c\p^{n} \]
is a one parameter family of H-minimal Lagrangian immersions.

We know from Lemma 2,3 when $\gamma_\theta $ is a closed curve, but now we want to
study when $\Pi \circ \phi_\theta$ is periodic
in its first variable. If we write $\gamma_\theta= (\rho_1 e^{i\nu_1}, \rho_2 e^{i\nu_2})$, Lemma 2,3 says that $\rho_i(t+T)=\rho_i (t)$, $i=1,2$.
Then it is not complicated to deduce that there exists $A>0$ such that $(\Pi \circ \phi_\theta)(t+A,p,q)=(\Pi \circ \phi_\theta)(t,p,q)$
if and only if there exists $\nu \in \r$  and $m\in \z$ ($A$ must be an integer multiple of $T$, $A=mT$) verifying
\begin{equation}\label{eq:per}
e^{i \nu_j (t+m T)} = e^{i\nu} e^{i \nu_j (t)}, \, j=1,2.
\end{equation}
From (\ref{eq:gamma0}) we can deduce that
\begin{equation}\label{eq:nu}
\rho_j^2 \nu_j' = (-1)^{j-1} c_\theta^{n_1+1} s_\theta^{n_2+1},\,  j=1,2.
\end{equation}
Then it is easy to check that $\nu_j(t+mT)=\nu_j(t)+m\nu_j (T)$, $j=1,2$, and (\ref{eq:per}) is equivalent to
$e^{i m \nu_j(T)}=e^{i\nu } $, $j=1,2$.
This means that $(\nu_2(T)-\nu_1(T))/2\pi$ must be a rational number. Using (\ref{eq:nu}) this implies that
\[
\theta \in \Gamma :=
\left\{ \theta \in (0,\frac{\pi}{2}) \, / \, \frac{\cos^{n_1+1}\theta \sin^{n_2+1}\theta}{2\pi} \int_0^{T}\frac{dt}{|\gamma_1|^2(t)|\gamma_2|^2(t)} \in\q \right\}.
\]
This study  leads to the following result.
\begin{corollary}
Given $\theta \in \Gamma$ and any  C-minimal Legendrian immersions
$\psi_i : N_i \longrightarrow \s^{2n_i+1}$, $i=1,2$, then the immersions $\phi_\theta$, $\theta \in \Gamma$, induce a one parameter family of H-minimal Lagrangian immersions
\[
\Phi_\theta : \s^1 \times N_1 \times N_2 \rightarrow \c\p^n, \, n=n_1+n_2+1, \, \theta \in \Gamma.
\]
In particular, $\Phi_\theta$ is minimal if and only if $\psi_i $, $i=1,2$, are minimal.
\end{corollary}

\vspace{1cm}

\section{H-minimal Lagrangian cones in complex \\ Euclidean space}

Let $\Omega_0=dz_1 \wedge \dots \wedge dz_{n+1} $ be the complex volume $(n\!+\!1)$-form on $\c^{n+1}$.
It is well-known that $\Re (e^{i\theta}\Omega_0)$, $\theta \in [0,2\pi )$, is the family of special Lagrangian calibrations in $\c^{n+1}$ (see [HL]) and that their calibrated submanifolds, the well known {\em special Lagrangian} submanifolds of $\c^{n+1}$, are not only minimal submanifolds but also minimizers in their homology class.

If $\psi:N^{n+1}\rightarrow \c^{n+1}$ is a Lagrangian immersion of an oriented manifold $N$ and
$\{ \widetilde{e}_0,\widetilde{e}_1,\dots,\widetilde{e}_n \}$ is an oriented orthonormal basis in $TN$, then the matrix
$\{ \psi_* (\widetilde{e}_0),\psi_*(\widetilde{e}_1),\dots,\psi_* (\widetilde{e}_n ) \}$ is an unitary matrix and so the following map is well defined:
\[
\begin{array}{c}
\widetilde{\beta} : N^{n+1} \longrightarrow \r / 2\pi \z
\\
e^{i\widetilde{\beta} (p)}=(\psi^* \Omega_0)_p (\widetilde{e}_0,\dots,\widetilde{e}_n) .
\end{array}
\]
$\widetilde{\beta}$ is known as the {\em Lagrangian angle} map of $\psi $ and verifies $J\widetilde{\nabla}\widetilde{\beta} = (n\!+\!1)\widetilde{H} $,
where $\widetilde{H}$ is the mean curvature of $\psi $. Then $\psi $ is a special Lagrangian immersion (with phase $\theta $)
if and only if $\widetilde{\beta}(p)=\theta$, $\forall p \in N$; moreover, $\psi $ is a Hamiltonian minimal Lagrangian immersion
if and only if $\widetilde{\beta}$ is a harmonic function.

Given a Legendrian immersion $\phi: M^n \rightarrow \s^{2n+1}$, the {\em cone} with link $\phi $ in $\c^{n+1}$ is the
map given by
\[
\begin{array}{c}
C(\phi): \r \times M^n \longrightarrow \c^{n+1}
\\
(s,p) \mapsto s \, \phi(p).
\end{array}
\]
It is clear that $C(\phi)$ is a Lagrangian immersion, i.e. $C(\phi)^*\omega \equiv 0$, with singularities at $s=0$.
We consider in what follows $s\neq 0$. The induced metric in $\r^* \times M$ by $C(\phi)$ is $ds^2 \times s^2 \langle, \rangle$, where $\langle, \rangle $ is the induced metric on $M$ by $\phi $. So, if $\{ e_1,\dots,e_n \}$ is an oriented orthonormal basis in $TM$, then $\{ \widetilde{e}_0=(1,0),\widetilde{e}_1=(0,\frac{e_1}{s}),\dots,\widetilde{e}_n=(0,\frac{e_n}{s}) \}$
is an oriented orthonormal frame on $T(\r^* \times M )$. Thus:
\[
\begin{array}{c}
e^{i \widetilde{\beta} (s,p)}=(\Omega_0)_{C(\phi)(s,p)} (C(\phi)_* (\widetilde{e}_0),\dots,C(\phi)_* (\widetilde{e}_n))=
\\ \\
\det_\c \{ \phi, \phi_* (e_1), \dots, \phi_* (e_n)  \}(p) = e^{i \beta (p)}.
\end{array}
\]
As a consequence, we deduce the following result.
\begin{proposition}
Let $\phi :M^n\rightarrow \s^{2n+1}$ be a Legendrian immersion of an oriented manifold $M$ and $C(\phi):\r \times M \rightarrow \c^{n+1}$ the cone with link $\phi$.
Then
$\phi $ is C-minimal if and only if $C(\phi)$ is H-minimal.
\end{proposition}
In particular, $\phi $ is minimal if and only if $C(\phi)$ is minimal; this result was used in [H1] and [H2].

Thanks to Proposition 3 we have a fruitful simple method of construction of examples of H-minimal Lagrangian cones in $\c^{n+1}$ using the C-mininal Legendrian immersions described in section 3. As an application of it, we finish this section with the following illustrative result.

\begin{corollary}
Let $\psi_i : N_i \longrightarrow \s^{2n_i+1}$, $i=1,2$, be C-minimal Legendrian immersions
of $n_i$-dimensional oriented Riemannian manifolds $N_i$, $i=1,2$, and $\gamma=(\gamma_1, \gamma_2): I \rightarrow \s^3 \subset \c^2$ a solution of some equation in the one parameter family of o.d.e.
\begin{equation}
\label{eq:gammamu2}
(\gamma_j' \overline{\gamma_j})(t) =
 (-1)^{j-1} i\, e^{i\mu t} \,    \overline{\gamma_1}(t)^{n_1+1} \,
\overline{\gamma_2}(t)^{n_2+1}, \,   \mu \in \r ,  \, j=1,2.
\end{equation}
Then
\[
\begin{array}{c}
\psi :\r \times I \times N_1 \times N_2 \rightarrow \c^{n+1}=\c^{n_1+1}\times\c^{n_2+1}, \,
  (n=n_1+n_2+1)
 \\ \\
\psi(s,t,p,q)= \left(
s\,\gamma_1 (t) \psi_1(p)  \, , \,
 s\,\gamma_2 (t) \psi_2(q)
\right)
\end{array}
\]
is a H-minimal Lagrangian cone.

In particular, if $\psi_i$, $i=1,2$, are minimal Legendrian immersions  and $\gamma$ is a solution of
(\ref{eq:gammamu2}) with $\mu =0$ (see Lemma 2), the
corresponding cone
$\psi: \r^2 \times N_1\times N_2 \longrightarrow \c^{n+1}$
is special Lagrangian.
\end{corollary}

A more explicit family of H-minimal Lagrangian cones is described when we consider the one parameter family of C-minimal Legendrian immersions coming from Corollary 2. Concretely,
for any $\delta \in (0,\pi/2)$ (denoting $c_{\delta}=\cos\delta, s_{\delta}=\sin\delta$),
\[
\begin{array}{c}
\psi_\delta: \r^2 \times N_1\times N_2 \longrightarrow \c^{n+1}=\c^{n_1+1}\times\c^{n_2+1} \,
\, (n=n_1+n_2+1)
\\  \\
\psi_\delta (s,t,p,q)=
( c_{\delta}\, s \, \exp (i s_{\delta}^{n_1+1} c_{\delta}^{n_2-1}   t) \,\psi_1 (p) \,  , \,
 s_{\delta}\, s \, \exp (-i s_{\delta}^{n_1-1} c_{\delta}^{n_2+1}   t) \,\psi_2 (q))
\end{array}
\]
is a H-minimal Lagrangian cone. In particular, if $\psi_i$, $i=1,2$, are minimal Legendrian immersions  and
$\tan^2 \delta_0 = (n_2+1)/(n_1+1)$, then $\psi_{\delta_0}$ is a special Lagrangian cone.

\section{The complex hyperbolic case}

In this section we summerize the analogous results when the
ambient space is  the complex hyperbolic space and we omit the
proofs of them.

Let $\c^{n+1}_1$ be the complex Euclidean space $\c^{n+1}$ endowed with the indefinite metric $\langle,\rangle =\Re\,(,)$, where
\[
(z,w)=\sum_{i=1}^n z_i\bar{w}_i-z_{n+1}\bar{w}_{n+1},
\]
for $z,w\in\c^{n+1}$, where $\bar z$ stands for the conjugate of $z$. The Liouville 1-form is given by $\Lambda_z (v) = \langle v, J z \rangle $, $\forall z \in \c^{n+1}, \forall v \in T_z \c^{n+1}$, and the Kaehler 2-form is $\omega = d\Lambda /2 $. We denote  by $\h^{2n+1}_1$ the anti-De Siter space, which is defined as the hypersurface of $\c^{n+1}_1$ given by
\[
\h^{2n+1}_1=\{z\in\c^{n+1}\,/\,(z,z)=-1\},
\]
and  by $\Pi:\h^{2n+1}\rightarrow\c\h^n$, $\Pi(z)=[z]$, the Hopf fibration of $\h^{2n+1}_1$ on the complex hyperbolic space $\c\h ^{n}$. We also denote the  metric, the complex structure and the K\"{a}hler two-form in $\c\h^n$ by $\langle, \rangle$, $J$ and $\omega$ respectively. This metric has constant holomorphic sectional curvature -4.
We will also note by $\Lambda $ the restriction to $\h^{2n+1}_1$ of the Lioville 1-form of $\c^{n+1}_1$.
So $\Lambda $ is the contact 1-form of the canonical (indefinite) Sasakian structure on the anti-De Sitter space $\h^{2n+1}_1$.
An immersion $\phi :M^n \rightarrow \h^{2n+1}_1$ of an $n$-dimensional manifold $M$ is said to be {\em Legendrian}
if $\phi^* \Lambda \equiv 0$. So $\phi $ is isotropic in $\c^{n+1}_1$, i.e. $\phi^* \omega \equiv 0$ and, in particular,
the normal bundle $T^\perp M = J (TM) \oplus {\rm span \, } \{ J\phi \} $. This means that $\phi $ is horizontal with respect to the Hopf fibration $\Pi:\h^{2n+1}_1 \rightarrow \c\h^n$ and, hence $\Phi= \Pi \circ \phi: M^n\rightarrow \c\h^n$ is a Lagrangian immersion and the
 induced metrics on $M^n$ by $\phi$ and $\Phi $ are the same.
It is easy to check that $J\phi $ is a totally geodesic normal vector field and so the second fundamental forms of $\phi $ and $\Phi $ are related by
\[
\Pi _* (\sigma_\phi (v,w)) = \sigma_\Phi (\Pi_* v, \Pi_* w) , \forall v,w \in TM .
\]
So the mean curvature vector $H$ of $\phi $ satisfies that $\langle H, J \phi \rangle =0$  and,
in particular, {\em $\phi :M^n \rightarrow \h^{2n+1}_1$ is minimal if and only if $\Phi= \Pi \circ \phi: M^n\rightarrow \c\h^n $ is minimal.}
In this way, we can construct (minimal) Lagrangian submanifolds in $\c\h^n$ by projecting Legendrian ones in $\h^{2n+1}_1$ by the Hopf fibration $\Pi$.
Let $\Omega $ be the complex $n$-form on $\h^{2n+1}_1$ given by
\[
\Omega_z(v_1,\dots,v_n)=\det_{\,\c} \, \{ z,v_1,\dots,v_n \} .
\]
If $\phi :M^n\rightarrow \h^{2n+1}_1$ is a Legendrian immersion of a manifold $M$,
then $\phi^* \Omega $ is a complex $n$-form on $M$. In the following result we analyze this $n$-form $\phi^* \Omega$.
\begin{lemma}
If $\phi :M^n\rightarrow \h^{2n+1}_1$ is a Legendrian immersion of a manifold $M$, then
\begin{equation} \label{eq:covariante2}
\nabla (\phi^* \Omega)= \alpha_H \otimes \phi^* \Omega,
\end{equation}
where $\alpha_H $ is the one-form on $M$ defined by $\alpha_H (v)= i\, n \langle H, Jv \rangle$
and $H$ is the mean curvature vector of $\phi$.
Consequently, if $\phi $ is minimal then $M$ is orientable.
\end{lemma}
\vspace{0.2cm}

Suppose that our Legendrian submanifold $M$ is oriented. Then we can consider the well defined map given by
\[
\begin{array}{c}
\beta : M^n \longrightarrow \r / 2\pi \z
\\
e^{i\beta (p)}=(\phi^* \Omega)_p (e_1,\dots,e_n)
\end{array}
\]
where $\{ e_1,\dots, e_n \}$ is an oriented orthonormal frame in $T_pM$. We will call $\beta$ the {\em Legendrian angle} map of $\phi $.
As a consequence of  (\ref{eq:covariante2})
we obtain
\begin{equation} \label{eq:gradiente2}
J\nabla \beta = n H,
\end{equation}
and so we deduce the following result.
\begin{proposition}
Let $\phi :M^n\rightarrow \h^{2n+1}_1$ be a Legendrian immersion of an oriented manifold $M$. Then $\phi $ is minimal if and only if the Legendrian angle map $\beta $ of $\phi $ is constant.
\end{proposition}
\vspace{0.3cm}

On the other hand, a vector field $X$ on $\h^{2n+1}_1$ is a {\em contact vector field} if ${\cal L}_X \Lambda = g \Lambda $, for some function
$g\in C^\infty (\h^{2n+1}_1)$, where $\cal L$ is the Lie derivative in $\h^{2n+1}_1$. Then $X$ is a contact vector field if and only if there exists $F\in C^\infty(\h^{2n+1}_1)$ such that
\[
X_z=J(\overline{\nabla}F)_z - 2F(z) Jz, \quad z\in \h^{2n+1}_1 ,
\]
where $\overline{\nabla}F$ is the gradient of $F$.
The diffeomorphisms of the flux $\{ \varphi_t \}$ of $X$ transform Legendrian submanifolds in Legendrian ones. In this setting, it is natural to study the following variational problem.
Let $\phi :M^n\rightarrow \h^{2n+1}_1$ a Legendrian immersion with mean curvature vector $H$. A normal vector field $\xi_f$ to $\phi $ is called  a {\em contact field} if
\[
\xi_f = J\nabla f - 2 f J\phi,
\]
where $f\in C^\infty (M)$ and $\nabla f$ is the gradient of $f$ respect to the induced metric.
If $f\in  C_0^\infty (M)$ and $\{ \phi_t :M\rightarrow\h_1^{2n+1} \}$ is a variation of $\phi $ with $\phi_0=\phi $ and
$\frac{d}{dt}_{\left|_{t=0}\right.}  \phi_t  = \xi_f $, the first variation of the volume functional is given by
\[
\frac{d}{dt}_{\left|_{t=0} \right.}  \, {\rm vol} (M,\phi_t^* \langle, \rangle )  = -\int_M  f\, {\rm div} \, JH \, dM .    \]
This means that the critical points of the above variational problem are Legendrian submanifolds such that
\[ {\rm div}  JH =0. \]
These critical points will be called {\em  contact minimal} (or briefly {\em C-minimal}) Legendrian submanifolds of $\h^{2n+1}_1$.
\begin{proposition}
If $\phi :M^n\rightarrow \h^{2n+1}_1$ is a Legendrian immersion of a Riemannian manifold $M$,
then:
\begin{enumerate}
\item If $M$ is oriented, $\phi $ is C-minimal if and only if the Legendrian angle $\beta $ of $\phi $ is a harmonic map.
\item $\phi $ is C-minimal if and only if $\Phi= \Pi \circ \phi :M^n\rightarrow \c\h^n$ is H-minimal.
\end{enumerate}
\end{proposition}
The identity component of the indefinite special orthogonal group will be denoted by $SO^1_0(m)$.
So $SO(n_1+1)\times SO_0^1(n_2+1)$ acts on $\h_1^{2n+1}\subset \c^{n+1}$, $n=n_1+n_2+1$, as a subgroup of isometries in the following way:
\begin{equation}
\label{eq:action2}
(A_1,A_2)\in SO(n_1+1)\times SO_0^1(n_2+1)\longmapsto \left(\begin{array}{c|c}
A_1 &  \\
\hline
&  A_2
\end{array}\right)\in SO^1_0(n+1).
\end{equation}
In the following we state (without proofs) the main results of section 3 adapted to this context.
\begin{theorem}\
Let $\psi_1 : N_1 \longrightarrow \s^{2n_1+1}\subset \c^{n_1+1}$ and $\psi_2 : N_2 \longrightarrow \h_1^{2n_2+1}\subset \c^{n_2+1}$ be Legendrian immersions
of $n_i$-dimensional oriented Riemannian manifolds $( N_i, g_{N_i} )$, $i=1,2$, and $\alpha=(\alpha_1,\alpha_2):I\rightarrow \h^3_1\subset \c^2$  be a Legendre curve.
The map
\[
\phi:I\times N_1 \times N_2 \longrightarrow \h_1^{2n+1}\subset \c^{n+1}=\c^{n_1+1}\times \c^{n_2+1}, \, n=n_1+n_2+1 ,
\]
 defined by
\[
\phi(s,p,q)=\left(\alpha_1(s)\psi_1(p),\alpha_2(s)\psi_2(q)\right)
\]
is a Legendrian immersion in $\h_1^{2n+1}$ whose induced metric is
\begin{equation}
\label{eq:metric2} \langle , \rangle = |\alpha'|^2 ds^2 +
|\alpha_1|^2 g_{N_1} + |\alpha_2|^2 g_{N_2}
\end{equation}
and whose Legendrian angle map is
\begin{equation}
\label{eq:angle2}
\beta_\phi \equiv n_1 \pi +\beta_\alpha + n_1 \arg \alpha_1 + n_2 \arg \alpha_2 + \beta_{\psi_1}
+\beta_{\psi_2} \quad \mod \, 2\pi,
\end{equation}
where $\beta_\alpha$ denotes the Legendre angle of $\alpha$ and $\beta_{\psi_i}$ the Legendrian angle map of $\psi_i$, $i=1,2$.

Moreover, a Legendrian immersion $\phi :M^n\longrightarrow \h_1^{2n+1}$  is  invariant under the action (\ref{eq:action2})
 of $SO(n_1+1)\times SO^1_0(n_2+1)$, with $n=n_1+n_2+1$ and $n_1,n_2 \geq 2$, if and only if
$\phi $ is locally congruent to one of the above Legendrian immersions when $\psi_1$ is the totally geodesic Legendrian embedding of $\,\s^{n_1}$ in $\s^{2n_1+1}$ and $\psi_2$ is the totally geodesic Legendrian embedding of $\,\r\h^{n_2}$ in $\h_1^{2n_2+1}$, where $\r\h^{n_2}=\{ (y_1,\dots,y_{n_2 +1}) \in \r^{n_2 +1} \, / \, \sum_{i=1}^{n_2}y_i^2 -y_{n_2+1}^2=-1, \, y_{n_2+1}>0 \}$ is the $n_2$-dimensional real hyperbolic space; that is,
$\phi $ is locally given by $\phi(s,x,y)=(\alpha_1(s)x,\alpha_2(s)y)$, $x\in \s^{n_1}$, $y\in \r\h^{n_2}$, for a certain Legendre curve $\alpha$ in $\h_1^3$.
\end{theorem}

\begin{remark}
{\rm
If $n_2=0$ (resp. $n_1=0$) in the above Theorem, projecting by the Hopf fibration $\Pi:\h_1^{2n+1} \rightarrow \c\h^n$,
we obtain the Examples 2 (resp. the Examples 3) given in [CMU1].
}
\end{remark}

\begin{corollary}\
Let $\psi_1 : N_1 \longrightarrow \s^{2n_1+1}\subset \c^{n_1+1}$ and $\psi_2 : N_2 \longrightarrow \h_1^{2n_2+1}\subset \c^{n_2+1}$ be C-minimal Legendrian immersions
of $n_i$-dimensional oriented Riemannian manifolds $N_i$, $i=1,2$, and $\alpha=(\alpha_1,\alpha_2):I\rightarrow \h^3_1\subset \c^2$  be a Legendre curve.
Then the Legendrian immersion described in Theorem 2 given by
\begin{eqnarray*}
\phi : I \times N_1\times N_2 \longrightarrow \h_1^{2n+1}, \, n=n_1+n_2+1,  \\
\phi(t,p,q)=\left(\alpha_1(t)\psi_1(p),\alpha_2(t)\psi_2(q)\right)
\end{eqnarray*}
is  C-minimal  if and only if, up to congruences,  $(\alpha_1,\alpha_2)$ is a solution of some equation in the one parameter family of o.d.e.
\begin{equation}
\label{eq:alpha12}
(\alpha_1' \overline{\alpha_1})(t) = ( \alpha_2' \overline{\alpha_2})(t)=
 i \, e^{i \mu t} \,    \overline{\alpha_1}(t)^{n_1+1} \,
\overline{\alpha_2}(t)^{n_2+1}, \,  \mu \in \r .
\end{equation}
Moreover, the above Legendrian immersion $\phi $ is minimal if and only if
$\psi_i$, $i=1,2$, are  minimal and $(\alpha_1,\alpha_2)$ is a solution of some o.d.e. of
(\ref{eq:alpha12}) with $\mu =0$.
\end{corollary}

If we consider the particular cases $n_2=0\Leftrightarrow n_1=n-1$ and $n_1=0\Leftrightarrow n_2=n-1$ in the minimal case of Corollary 7, we recover (projecting via the Hopf fibration $\Pi$) the minimal Lagrangian submanifolds of $\c\h^n$ described in [CMU2, Propositions 3 and 5], although we used there an unit speed parametrization for $\alpha$.

From these two last results we can get similar examples to the ones given in Section 4 in the projective case.
Concretely, it is easy to check that for any $\rho >0$ the Legendre curve
\begin{equation}\label{eq:alfaro}
\alpha_\rho (t)=( sh_{\rho} \, \exp(i\, sh_{\rho}^{n_1-1}  ch_{\rho}^{n_2+1}    t) ,
ch_{\rho} \, \exp (i \, sh_{\rho}^{n_1+1} ch_{\rho}^{n_2-1}  t)),
\end{equation}
satisfies (\ref{eq:alpha12}) for
$\mu = sh_{\rho}^{n_1-1} ch_{\rho}^{n_2-1} \left( (n_1+1)ch_{\rho}^2  + (n_2+1)sh_{\rho}^2  \right)$, where $ch_{\rho}=\cosh\rho$, $sh_{\rho}=\sinh\rho$.

Hence an analogous reasoning like in Corollary 3 let us to obtain the following explicit family of examples.

\begin{corollary}\
Let $\psi_1 : N_1 \longrightarrow \s^{2n_1+1}\subset \c^{n_1+1}$ and $\psi_2 : N_2 \longrightarrow \h_1^{2n_2+1}\subset \c^{n_2+1}$ be C-minimal Legendrian immersions
of $n_i$-dimensional Riemannian manifolds $N_i $, $i=1,2$, and $\rho >0$. Then
\[
\begin{array}{c}
\Phi_\rho: \s^1 \times N_1\times N_2 \longrightarrow \c\h^{n}, \,
\, n=n_1+n_2+1,
\\  \\
\Phi_\rho (e^{is},p,q)= [( \sinh \rho \, \exp(i s\, \cosh^2\rho )
\,\psi_1 (p) \,  , \, \cosh \rho \, \exp (i s\, \sinh^2 \rho)
\,\psi_2 (q))]
\end{array}
\]
is a H-minimal Lagrangian immersion.
\end{corollary}
A particular case of Corollary 8 gives a one parameter family of H-minimal Lagrangian embeddings.

\begin{corollary}
For each $\rho >0$, the immersion $\Phi_{\rho}$ (given in Corollary 8) where $\psi_1$ (resp. $\psi_2$) is the totally geodesic Legendrian embedding of \,$\s^{n_1}$ into $\s^{2n_1+1}$ (resp. of \,$\r\h^{n_2}$ into $\h_1^{2n_2+1}$), provides a H-minimal Lagrangian embedding
\[
\begin{array}{c}
\frac{\textstyle \s^1\times\s^{n_1}\times\r\h^{n_2}}{\textstyle \z_2} \longrightarrow \c\h^n,\, n=n_1+n_2+1 \\ \\
\overline{(e^{is},x,y)} \longmapsto  [( \sinh \rho \, \exp(i s\,
\cosh^2 \rho ) \,x \,  , \, \cosh \rho \, \exp (i s\, \sinh^2
\rho)  \,y)]
\end{array}
\]
of the quotient of $\s^1\times\s^{n_1}\times\r\h^{n_2}$ by the action of the group $\z_2$
where the generator $h$ of $\z_2$ acts on $\s^1\times\s^{n_1}\times\r\h^{n_2}$ in the following way
\[
h(e^{is},x,y)=(-e^{is},-x,y).
\]
\end{corollary}

\vspace{0.1cm}

We finally pay now our attention to the equation (\ref{eq:alpha12}) with $\mu=0$.
We observe that it is exactly equation (3) in
[CU2, Lemma 2] (in the notation of that paper, put $p=n_1$ and $q=n_2$).
If we choose the initial conditions $\alpha(0)=(\sinh \varrho, \cosh \varrho)$, $\varrho >0$, we can make use of the study made in [CU2].

\begin{lemma}
Let $\alpha_\varrho=(\alpha_1,\alpha_2):I\subset \r \rightarrow \h_1^3$ be the only curve solution of
\begin{equation}\label{eq:alfa0}
\alpha'_j\bar{\alpha}_j= i\, \bar{\alpha}_1^{n_1+1}\bar{\alpha}_2^{n_2+1}, \, j=1,2,
\end{equation}
satisfying the real initial conditions
$\alpha_\varrho(0)=(\sinh \varrho, \cosh \varrho)$, $\varrho >0$.
Then:
\begin{enumerate}
\item $\Re (\alpha_1^{n_1+1}\alpha_2^{n_2+1})=\sinh^{n_1+1} \varrho \cosh^{n_2+1} \varrho$.
\item For $j=1,2$,
$\bar{\alpha}_j(t)=\alpha_j(-t),\,\forall t\in I$.
\item The curves $\alpha_j,\,j=1,2$, are embedded and can be parameterized by $\alpha_j(s)=\rho_j(s)e^{i\theta_j(s)}$, where
\begin{eqnarray*}
\rho_1(s)&=&\sqrt{s^2+sh_\varrho^2 },\\
 \theta_1(s)&=&\int_0^s\frac{sh_\varrho^{n_1+1} ch_\varrho^{n_2+1} x\,dx}{(x^2+sh_\varrho^2 )
\sqrt{(x^2+sh_\varrho^2 )^{n_1+1}(x^2+ch_\varrho^2 )^{n_2+1}-sh_\varrho^{2(n_1+1)} ch_\varrho^{2(n_2+1)}}}
\end{eqnarray*}
and
\begin{eqnarray*}
\rho_2(s)&=&\sqrt{s^2+ch_\varrho^2},\\
 \theta_2(s)&=&\int_0^s\frac{sh_\varrho^{n_1+1} ch_\varrho^{n_2+1} x\,dx}{(x^2+ch_\varrho^2 )
\sqrt{(x^2+sh_\varrho^2 )^{n_1+1}(x^2+ch_\varrho^2 )^{n_2+1}-sh_\varrho^{2(n_1+1)} ch_\varrho^{2(n_2+1)}}},
\end{eqnarray*}
where $ch_{\varrho}=\cosh\varrho$, $sh_{\varrho}=\sinh\varrho$.
\end{enumerate}

\end{lemma}

 In this way,  the immersions $\phi_\varrho$, $\varrho >0$, constructed with the curves $\alpha_\varrho$ of Lemma 4 induce a one parameter family of H-minimal Lagrangian immersions
\[
\Phi_\varrho : \r \times N_1 \times N_2 \rightarrow \c\h^n, \, n=n_1+n_2+1, \, \varrho >0.
\]
In particular, $\Phi_\varrho$ is minimal if and only if $\psi_i $, $i=1,2$, are minimal.
We conclude with the following particular case that leads to a one parameter family of minimal Lagrangian embeddings.
\begin{corollary}
Let $\varrho >0$ and denote $ch_{\varrho}=\cosh\varrho$, $sh_{\varrho}=\sinh\varrho$.
Then
\[
\begin{array}{c}
\r \times \s^{n_1} \times \r\h^{n_2} \longrightarrow \c\h^{n}, \,
\, n=n_1+n_2+1,
\\  \\
(s,x,y) \mapsto
[( \sqrt{s^2+sh_\varrho^2 } \, \exp (i\, \theta_1 (s)) x \,  , \,
\sqrt{s^2+ch_\varrho^2}   \, \exp (i\, \theta_2 (s)) y)],
\end{array}
\]
where $\theta_i(s)$, $i=1,2$, are given in part 3 of Lemma 4, is a
minimal Lagrangian embedding.
\end{corollary}

\vspace{1cm}

{\sc addresses}:

(first author)

Departamento de Matem\'{a}ticas

Escuela Polit\'{e}cnica Superior

Universidad de Ja\'{e}n

23071 Ja\'{e}n

SPAIN

{\tt icastro@ujaen.es}

(second author)

Department of Mathematical Sciences

Tsinghua University

100084 Beijing

PEOPLE'S REPUBLIC OF CHINA

{\tt hli@math.tsinghua.edu.cn}

(third author)

Departamento de Geometr\'{\i}a  y Topolog\'{\i}a

Universidad de Granada

18071 Granada

SPAIN

{\tt furbano@ugr.es}

\end{document}